\documentclass{article}
\usepackage[notcite, notref,color,final]{showkeys}
\definecolor{refkey}{rgb}{0,0,1}
\definecolor{labelkey}{rgb}{0,0,1}
\usepackage{amsmath,amssymb,latexsym}
\usepackage[amsmath,thmmarks,thref]{ntheorem}
\usepackage{xspace,amsfonts,latexsym}
\usepackage{graphicx}
\usepackage[numbers]{natbib}
\usepackage{hyperref}
\include{thmstuff}
\newcommand{\boldp}{\ensuremath{{\bf p}}\xspace}
\newcommand{\boldq}{\ensuremath{{\bf q}}\xspace}
\newcommand{\boldpq}{\ensuremath{{\bf p};\;{\bf q}}\xspace}
\newcommand{\boldptimesq}{\ensuremath{{\bf p \times q}}\xspace}
\newcommand{\boldpminq}{\ensuremath{{\bf p; -q}}\xspace}
\newcommand{\chihat}{\ensuremath{\widehat{\chi}}\xspace}

\newcommand{\boldptq}{\ensuremath{ {\bf p \times q}}\xspace}
\newcommand{\colourS}{\ensuremath{\mfS^{(m)}}\xspace}
\newcommand{\km}{\ensuremath{\kappa^{(m)}}\xspace}
\newcommand{\gent}{{\ensuremath{T_{\boldpq}(x)}}\xspace}
\newcommand{\geni}{{\ensuremath{I_{\boldpq}(x)}}\xspace}

\newcommand{\tclass}{\ensuremath{\mathcal{T}}\xspace}
\newcommand{\bclass}{\ensuremath{\mathcal{B}}\xspace}
\newcommand{\wclass}{\ensuremath{\mathcal{W}}\xspace}
\newcommand{\wclasshat}{\ensuremath{\hat{\mathcal{W}}}\xspace}
\newcommand{\what}{\ensuremath{\hat{W}}\xspace}
\newcommand{\fact}{\ensuremath{\mathrm{TopFact}}\xspace}
\begin{document}
\title{Stanley's character polynomials and coloured factorizations in the
symmetric group}
\author{A. Rattan\\Department of Mathematics\\Massachusetts
Institute of Technology\\Cambridge, MA, 02139\\ email: {\small\texttt
arattan@math.mit.edu}}
\date{February 10, 2007}

\newcommand{\seqtwo}[2]{#1_1, \ldots, #1_{#2}}
\newcommand{\seqinf}[1]{#1_1, #1_{2}, \ldots}
\newcommand{\seq}[2]{\ensuremath{#1_1, #1_2, \ldots, #1_{#2}}}
\newcommand{\seqo}[2]{#1_0, #1_1, \ldots, #1_{#2}}
\newcommand{\seqtwofromto}[3]{#1_{#2}, \ldots, #1_{#3}}
\newcommand{\seqfromto}[3]{#1_{#2}, #1_{#2 + 1}, \ldots, #1_{#3}}
\newcommand{\thbox}{\hfill$\Box$}
\newcommand{\parti}[1]{1^{#1_1}2^{#1_2} \ldots n^{#1_n}}
\newcommand{\parts}[1]{m_1(#1), m_2(#1), \ldots, m_n(#1)}
\newcommand{\partso}[1]{m_0(#1), m_1(#1), \ldots, m_n(#1)}
\newcommand{\set}[2]{\{\, #1 \; |\; #2\}}
\newcommand{\diffpart}[3]{#1_1^{#2_1}#1_2^{#2_2} \ldots #1_{#3}^{#2_{#3}}}
\newcommand{\sumofvars}[2]{#1_1 + #1_2 + \cdots + #1_{#2}}
\newcommand{\sumofvarstwo}[2]{#1_1 + \cdots + #1_{#2}}
\newcommand{\diffofvars}[2]{-#1_1 - #1_2 - \cdots - #1_{#2}}
\newcommand{\lessvars}[2]{#1_1 < #1_2 < \cdots < #1_{#2}}
\newcommand{\greatvars}[2]{#1_1 > #1_2 > \cdots > #1_{#2}}
\newcommand{\prodofvars}[2]{#1_1  #1_2 \cdots #1_{#2}}
\newcommand{\prodtwofromto}[3]{#1_{#2} \cdots #1_{#3}}
\newcommand{\upst}{\mathrm{st}}
\newcommand{\upnd}{\mathrm{nd}}
\newcommand{\upth}{\mathrm{th}}
\newcommand{\coeff}[1]{\ensuremath{[#1]\;}\xspace}

\newcommand{\seqpermvars}[3]{#1_{#2(1)}, #1_{#2(2)}, \ldots, #1_{#2(#3)}}
\newcommand{\smalldisunion}{\ensuremath{\stackrel{\cdot}{\cup}}}
\newcommand{\bigdisunion}{\ensuremath{\stackrel{\cdot}{\bigcup}}}
\newcommand{\symgroup}[1]{\ensuremath{\mathfrak{S}_#1}}
\newcommand{\fullcyc}[1]{\ensuremath{(1 \, 2 \, \cdots \, #1)}}
\newcommand{\laginv}{ {\ensuremath{\langle -1 \rangle}}}
\newcommand{\onenk}{1^{n-k}} 
\newcommand{\supp}{\ensuremath{\mathrm{supp}}\xspace}
\newcommand{\stm}{\ensuremath{\setminus}\xspace}
\newcommand{\al}{\ensuremath{\alpha}\xspace}
\newcommand{\be}{\ensuremath{\beta}\xspace}
\newcommand{\ga}{\ensuremath{\gamma}\xspace}
\newcommand{\lam}{\ensuremath{\lambda}\xspace}
\newcommand{\si}{\ensuremath{\sigma}\xspace}
\newcommand{\om}{{\ensuremath{\omega}\xspace}}
\newcommand{\vt}{\ensuremath{\vartheta}\xspace}
\newcommand{\sms}{\ensuremath{\setminus}\xspace}
\newcommand{\rar}{\ensuremath{\rightarrow}\xspace}
\newcommand{\mpt}{\ensuremath{\mapsto}\xspace}
\newcommand{\f}{\ensuremath{\frac}}
\newcommand{\tf}{\ensuremath{\tfrac}}
\newcommand{\ld}{\ensuremath{\ldots}\xspace}
\newcommand{\mbC}{\ensuremath{\mathbb{C}}\xspace}
\newcommand{\mbR}{\ensuremath{\mathbb{R}}\xspace}
\newcommand{\mcC}{\ensuremath{\mathcal{C}}\xspace}
\newcommand{\mcS}{\ensuremath{\mathcal{S}}\xspace}
\newcommand{\mcR}{\ensuremath{\mathcal{R}}\xspace}
\newcommand{\mcT}{\ensuremath{\mathcal{T}}\xspace}
\newcommand{\mcB}{\ensuremath{\mathcal{B}}\xspace}
\newcommand{\mcP}{\ensuremath{\mathcal{P}}\xspace}
\newcommand{\mfS}{\ensuremath{\mathfrak{S}}\xspace}


\newcommand{\grob}{Gr\"{o}bner basis }
\newcommand{\grobs}{Gr\"{o}bner bases }
\newcommand{\ky}[1]{k[y_1, y_2, \ldots, y_{#1}]}
\newcommand{\kx}[1]{k[x_1, x_2, \ldots, x_{#1}]}
\newcommand{\kyx}[2]{k[y_1, \ldots, y_{#1}, x_1, \ldots, x_{#2}]}
\newcommand{\kxy}[2]{k[x_1, \ldots, x_{#1}, y_1, \ldots, y_{#2}]}
\newcommand{\id}[1]{\langle #1 \rangle}
\newcommand{\qedend}{\hfill\qed}
\newcommand{\isubnk}{I_n^{(k)}}
\newcommand{\psubnk}{P_n^{(k)}}
\newcommand{\khat}{\hat{k}}
\newcommand{\kpark}{{(k)}}
\newcommand{\Phat}{\hat{P}}
\newcommand{\Ihat}{\hat{I}}


\maketitle

\begin{abstract}
	\noindent In \citet{stan:7}, the author introduces polynomials which help
	evaluate symmetric group characters and conjectures
	that the coefficients of the polynomials are positive.  In
	\cite{stan:9},
	the same author gives a conjectured combinatorial interpretation for the
	coefficients of the polynomials.  Here, we prove the conjecture for the
	terms of highest degree.
\end{abstract}

\section{Introduction}

A {\em partition} is a weakly ordered list of positive
integers $\lam =\lam_1\lam_2\ld\lam_k$, where $\lam_1\geq\lam_2\geq\ld\geq\lam_k$.
The integers $\lam_1,\ld ,\lam_k$ are called the {\em parts} of
the partition $\lam$, and we denote the number of parts by $\ell(\lam)=k$.
If $\lam_1+\ld +\lam_k=d$, then $\lam$ is a partition of $d$, and
we write $\lam\vdash d$. We denote by $\mcP$ the set of all
partitions, including the single partition of $0$ (which
has no  parts).  For partitions $\om, \lam \vdash n$, let $\chi_\om(\lam)$ be
the character of the irreducible representation of the symmetric
group $\mfS_n$ indexed by $\om$, and evaluated on the conjugacy class
 $\mcC_{\lam}$ of $\mfS_n$, where $C_\lambda$ is the class of all permutations
 whose disjoint cycle lengths are specified by the parts of $\lam$.  For
 a permutation $\alpha$, we use the notation $\kappa(\alpha)$ to denote the number of
 cycles of $\alpha$.  We use the convention that permutation are multiplied from
 right to left.

Various scalings of irreducible symmetric group characters have been considered in
the recent literature, one of which is the central object of this paper.  
Suppose that $\mu \vdash k$ and $k\leq n$.  For the conjugacy class $\mcC_{\mu
1^{n-k}}$, the {\em normalized character} is given by
\begin{equation*}
  \chihat_\om(\mu\; 1^{n-k}) = (n)_k
\frac{ \chi_\om(\mu\; 1^{n-k})}{\chi_\om(1^n)},
\end{equation*}
where $(n)_k$ is the falling factorial $n (n-1) \cdots (n-k+1)$.
The normalized character has been the topic of much recent literature and has
been shown to have connections with combinatorics and free probability, see for
example \cite{bi, bi1, rattgoul:3, sniasymp, stan:7}.

\begin{figure}[ht]
 \centering
  \includegraphics[width=6.0cm]{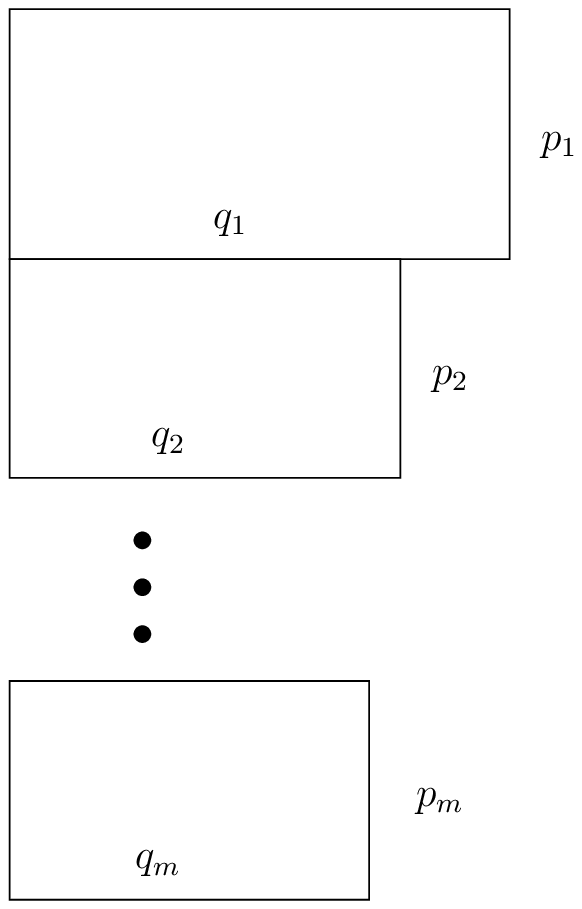}
  \caption{The shape \boldptimesq.}\label{fig:genpart}
\end{figure}

The subject of this paper is a particular polynomial expression for the
normalized character, introduced in \citet{stan:7}.  Consider the 
partition of $n$ with $p_i$ parts of size $q_i$, for $i$ from 1 to $m$, with $q_1$ the
largest part.
Thus, $\seq{p}{m}$ are
positive integers and $q_1 \geq q_2 \geq \cdots \geq q_m$ (see Figure
\ref{fig:genpart}).  We denote this partition of $n$ by \boldptimesq.  
Define the expression $F_k$ in indeterminates $\seqtwo{p}{m}, \seqtwo{q}{m}$ by
\begin{equation}\label{eq:deffk}
	F_k(\seq{p}{m};\; \seq{q}{m}) = \chihat_\boldptimesq(k\; 1^{n-k}).
\end{equation}
We often use $\boldp$ for 
$(\seqtwo{p}{m})$ and $\boldq$ for $(\seqtwo{q}{m})$, giving us the
notation $F_k(\boldp;\; \boldq)$
for $F_k(\seq{p}{m};\; \seq{q}{m})$.  The following theorem appears in
\citet[Proposition 1]{stan:7}.
\begin{theorem}[Stanley]
$F_k(\boldpq)$ is a polynomial in the $p$'s and $q$'s such that
\linebreak $(-1)^k F_k(1,1,\dots,1;\; -1,-1,\dots,-1) = (k+m-1)_k$. 
\end{theorem}
In light of this theorem, we call the polynomials in \eqref{eq:deffk}
\emph{Stanley's character polynomials}.  These polynomials are the main objects
in this paper.  For example, for the case $m=2$, the first two 
polynomials are
\begin{align*}
F_1(a,p;\; b, q) &= -ab-pq,\\
F_2(a,p;\; b, q) &= -{a}^{2}b+a{b}^{2}-2\,apq-{p}^{2}q+p{q}^{2}
\end{align*}
where we have set $p_1 = a, p_2 = p, q_1 = b$ and $q_2 =q$.
Also in \citet{stan:7} the author states that if one defines $F_\mu(\boldp;\;
\boldq)$ as
\begin{equation}\label{eq:gendef}
	F_\mu(\boldp;\; \boldq) = \chihat_\boldptimesq(\mu\; 1^{n-k}),
\end{equation}
where $\mu \vdash k \leq n$ then $F_\mu(\boldp;\; \boldq)$ is also polynomial.
We emphasize that $F_\mu(\boldpq)$ is independent of $n$ as it is, formally, a
polynomial in indeterminates in $\seqtwo{p}{m}$ and $\seqtwo{q}{m}$.  Our
understanding is that the polynomial $F_\mu(\boldpq)$ evaluates to the normalized
character in \eqref{eq:gendef} when evaluated at a shape $\boldptq$ that is a
partition of $n \geq k$.

In \citet{stan:7}, the author conjectures that $(-1)^k F_\mu(\boldpminq)$
has positive coefficients.  This conjecture has only been proved in the
case $m=1$ for general $\mu$ (see \cite[Theorem 1.1]{stan:7});  in particular, the conjecture is not known to
be true for $m > 1$ even when $\mu$ has one part; that is, the conjecture is unknown
even for $(-1)^k F_k(\boldpminq)$.  Some partial results showing positivity of
the coefficients of $(-1)^k F_k(\boldpminq)$ were given in
\citet{ratt:6}, but otherwise little is known about these polynomials.
Recently, Stanley \cite{stan:9} has a conjectured combinatorial interpretation for
$F_\mu(\boldp;\; \boldq)$, which we now explain.

Let $[m]$ be the set $\{1, 2, \dots, m\}$ and $\colourS_k$ be the
set of permutations of the set $[k]$ whose cycles are coloured by $[m]$.
Formally, if $C_k(\alpha)$ is the set of cycles in $\alpha$ then members of
$\colourS_k$ are ordered pairs $(\alpha, \psi)$ where $\alpha \in \mfS_k$ and
$\psi : C_k(\alpha) \longrightarrow [m]$.  Define a product $\circ : \colourS_k
\times \mfS_k \longrightarrow \colourS_k$ by the following:  for $(\alpha, \psi)
\in \colourS_k$, $\beta \in \mfS_k$ and $(\alpha, \psi) \circ \beta = (\gamma,
\nu)$, where
\begin{enumerate}
	\item $\gamma = \alpha \beta$, and
	\item If $u = (u_1\; u_2\; \cdots u_t)$ is a cycle of $\gamma$ and
		$C^\alpha(u_i)$ is the cycle of $\alpha$ containing the symbol
		$u_i$ then
		\begin{equation*}
			\nu(u) = \max_{1 \leq i \leq t} \{\psi(C^\alpha(u_i))\}.
		\end{equation*}
\end{enumerate}
In words, $\nu(u) = t$, where $t$ is the largest value of $\psi(w)$, and where $w$
ranges over all cycles in $\alpha$ with an element in common with $u$ (see
\cite[Page 3]{stan:9} for an example).
For $(\alpha, \psi) \in \colourS_k$ let $\km(\alpha, \psi) =
(\km_1(\alpha, \psi), \km_2(\alpha, \psi), \ldots)$ where $\km_i(\alpha, \psi)$
is the number of cycles of $\alpha$ coloured $i$ and for $\boldp =
(\seq{p}{m})$, $\boldp^{\km(\alpha, \psi)} = \prod_i p_i^{\km_i(\alpha, \psi)}$.
We can now state the conjecture we address in this paper, found in
\citet{stan:9}.
\begin{conjecture}[Stanley]\label{thm:mainstan}
	Suppose that $\mu \vdash k$, and let $\omega_\mu$ be some fixed element in the conjugacy class $C_\mu$ in $\mfS_k$.  Then
	\begin{equation*}
		(-1)^k F_\mu({\bf p};\;{\bf -q}) = \sum_{(\alpha, \psi) \in \colourS_k}
		\boldp^{\km(\alpha, \psi)} \boldq^{\km( (\alpha, \psi)\circ
		\omega_\mu)}.
	\end{equation*}
\end{conjecture}
As stated earlier, in \citet[Theorem 1.1] {stan:7} Conjecture \ref{thm:mainstan} has been proved
for $m=1$ (note that this corresponds to factorizations without any colours),
but otherwise this conjecture remains open.  However, for arbitrary
$m$, it is shown in \cite{stan:7} that
\begin{equation}\label{eq:fsubk}
	F_k(\boldp;\; \boldq)= -\f{1}{k} [x^{-1}]_\infty (x)_k \; \prod_{j=1}^m\f{
	\;  \left(x - (q_j + p_j + p_{j+1} + \cdots +
	p_m)\right)_k}{ \left(x - (q_j + p_{j+1} +
	p_{j+2} + \cdots + p_m)\right)_k},
\end{equation}
where for an expression $g(x)$ the notation $[x^{-1}]_\infty g(x)$ is the coefficient of $1/x$ when $g(x)$ is
expanded in powers of $1/x$.  From this it follows (see \cite[Proposition
2]{stan:7}) that if $G_k(\boldpq)$ are the terms of highest degree in
$F_k(\boldpq)$ and $G_{\boldpq}(x) = 1 + \sum_{k=0}^\infty G_k(\boldpq) x^{k+1}$
then
\begin{equation}
	G_{\boldpq}(x) = 1 + \sum_{i \geq 1} G_{i-1}(\boldp;\; \boldq)
	x^i 
	= \f{x}{\left(x\prod_{j=1}^m\f{
	\left(1-\left(q_j + 
	p_{j+1} + \cdots + p_m \right)x\right)}{\left(1 - \left(q_j +
	 p_j + \cdots + p_m\right)x\right)} \right)^\laginv}\label{eq:stantop}.
\end{equation}
where $\laginv$ denotes compositional inverse.  It easily follows from
\eqref{eq:stantop} that $G_{\boldpq}(x)$ satisfies
\begin{equation}\label{eq:stantopgoodform}
	-G_{\boldpminq}(-x) \prod_{j=1}^m \frac{\left(-G_{\boldpminq}(-x) -
	(p_j + \cdots + p_m)x + q_jx\right)}{\left( -G_{\boldpminq}(-x) -
	(p_{j+1} + \cdots + p_m)x + q_j x\right)} = -1
\end{equation}
It is \eqref{eq:stantopgoodform} that we will eventually use to prove our main
theorem.  It is known for $m=1$ (\citet[Page 9]{stan:7}) that the
series $G_{p;\; q}(x)$ is the generating series for \emph{top} factorizations in the symmetric group;  namely, we have
\begin{equation*}
	G_{p;\; q}(x) = 1 + px + \sum_{k \geq 1} x^{k+1} \sum_{{u \in \mfS_k \atop
	\kappa(u) + \kappa(u \omega_k) = k+1}} (-1)^k p^{\kappa(u)}(-q)^{\kappa(u \cdot
	\omega_k)}.
\end{equation*}
where $\omega_k$ is used for $\omega_{(1\; 2\; \cdots\; k)}$.
Here, we have $\kappa(u) = \kappa^{(1)}(u)$ is the number of cycles of $u$.
We call products of the type in the previous sum i.e. products of permutations
$\alpha \beta = \gamma$ in $\mfS_k$ such $\kappa(\alpha) + \kappa(\beta) = k
 + \kappa(\gamma)$ \emph{top products}, \emph{top factorizations} or \emph{minimal
factorizations}.  Such factorizations are an extremal case;  namely, if
$\alpha, \beta, \gamma \in \mfS_k$ and $\alpha \beta = \gamma$ then
\begin{equation}\label{eq:minfacts}
	\kappa(\alpha)+ \kappa(\beta) \leq k + \kappa(\gamma)
\end{equation}
(see \citet{goul:7}).  

Set,
\begin{equation*}
	\fact_{\boldpq}(x) = \sum_{k \geq 1} x^{k+1}\sum_{(\alpha, \psi) \in
	\colourS_k \atop \kappa(\alpha) + \kappa(\alpha \omega_k) = k+1} 
	\boldp^{\km(\alpha, \psi)}\boldq^{\kappa^{(m)}( (\alpha, \psi) \circ
	\omega_k)}.
\end{equation*}
The following are the two main theorems of this paper.  We will see that
Corollary \ref{thm:main2} follows from Theorem \ref{thm:main1}.
\begin{theorem}[Main Theorem]\label{thm:main1}
	Conjecture \ref{thm:mainstan} holds for the term of highest 
	degree in $(-1)^kF_k(\boldp;\; -\boldq)$;  that is,
\begin{equation}\label{eq:gentopterms}
	-G_{\boldpminq}(-x) = -1 + (\sumofvars{p}{m})x + \fact_{\boldpq}(x).
\end{equation}
\end{theorem}

\begin{cor}[Main Corollary]\label{thm:main2}
	For any partition $\mu \vdash k$, Conjecture \ref{thm:mainstan} holds for the
	terms of highest degree in $(-1)^kF_\mu(\boldp;\; -\boldq)$.
\end{cor}
We prove the main theorems at the end of Section \ref{sec:proofmain}.

\section{The Goulden-Jackson construction for top
factorizations}\label{sec:goulcon}

\begin{figure}[ht]
 \centering
  \includegraphics[width=14.0cm]{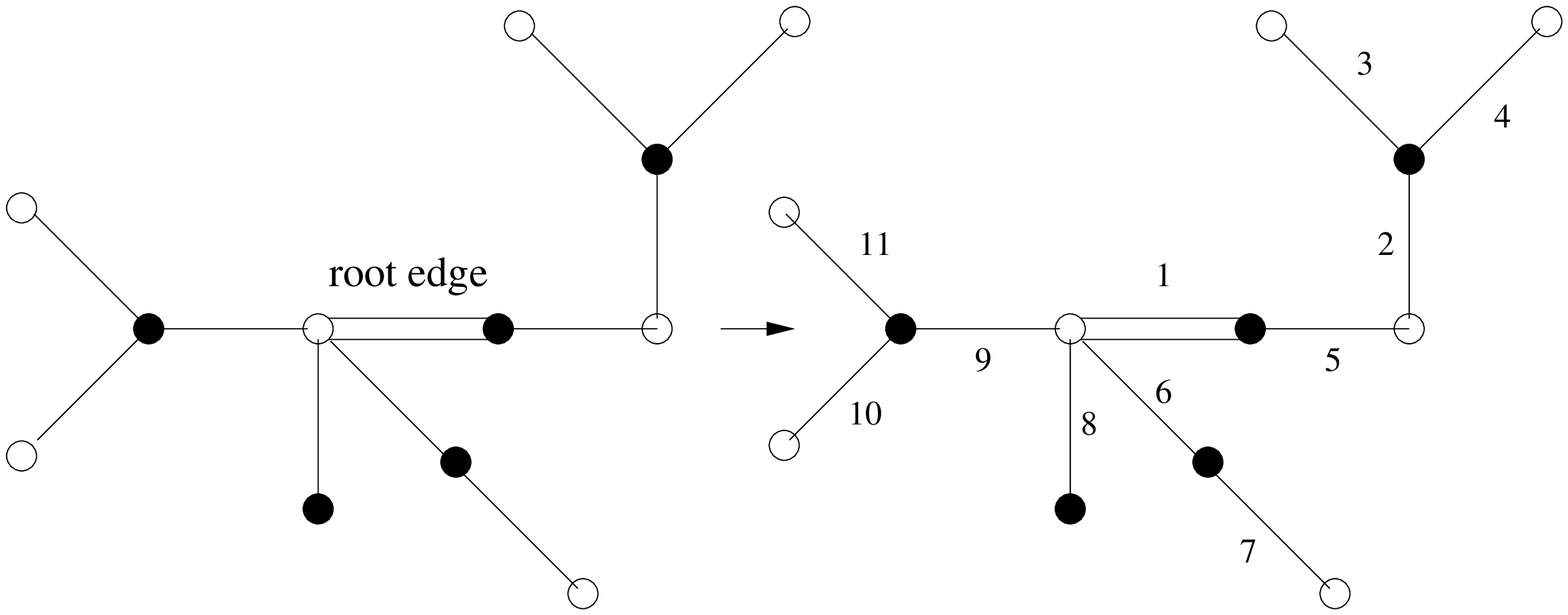}
  \caption{On the left is a black and white plane edge rooted tree.  Using the
  description in the first paragraph of Section \ref{sec:goulcon} to label the edges
  of the tree on the left, we obtain
  the tree on the right.  A clockwise rotation around each white vertex gives a cycle
  of the permutation $\alpha =
  (1\; 6\; 8\; 9)(2\; 5)(3)(4)(7)(10)(11)$, and likewise for the black vertices and $\beta = (1\; 5)(2\; 3\; 4)(6\;
  7)(8)(9\; 10\; 11)$.}\label{fig:baseplanetree}
\end{figure}

In \citet{goul:7}, the authors give a construction for top factorizations in the
symmetric group in terms of black and white plane edge rooted trees.  Namely, they
give a bijection between products of permutations $\alpha \beta = (1\; 2\;
\cdots k)$ in $\mfS_k$ such that $\kappa(\alpha) + \kappa(\beta) = k +1$ and edge
rooted plane 
trees on $k+1$ vertices, with vertices coloured black and white such that adjacent vertices
receive different colours.  The correspondence is very simple to state;  in the
tree, label the edges beginning with the root edge (which obtains the label 1).  The
edges are labelled in numerical order by travelling around the tree, keeping the
tree to the right and labelling an edge only when traversed from its white vertex to its
black vertex.  From each white vertex a cycle is obtained for the permutation
$\alpha$ by considering the sequence of edges incident with the vertex in a
clockwise direction.  Likewise, the cycles of $\beta$ are obtained from the black
vertices.  In Figure \ref{fig:baseplanetree}, the plane tree given corresponds
to the pair $\alpha = (1\; 6\; 8\; 9)(2\; 5)(3)(4)(7)(10)(11)$ and $\beta =
(1\; 5)(2\; 3\; 4)(6\; 7)(8)(9\; 10\; 11)$.  One can easily check that
$\alpha \beta = (1\; 2\; \cdots\; 11)$.

\begin{figure}[ht]
 \centering
  \includegraphics[width=8.0cm]{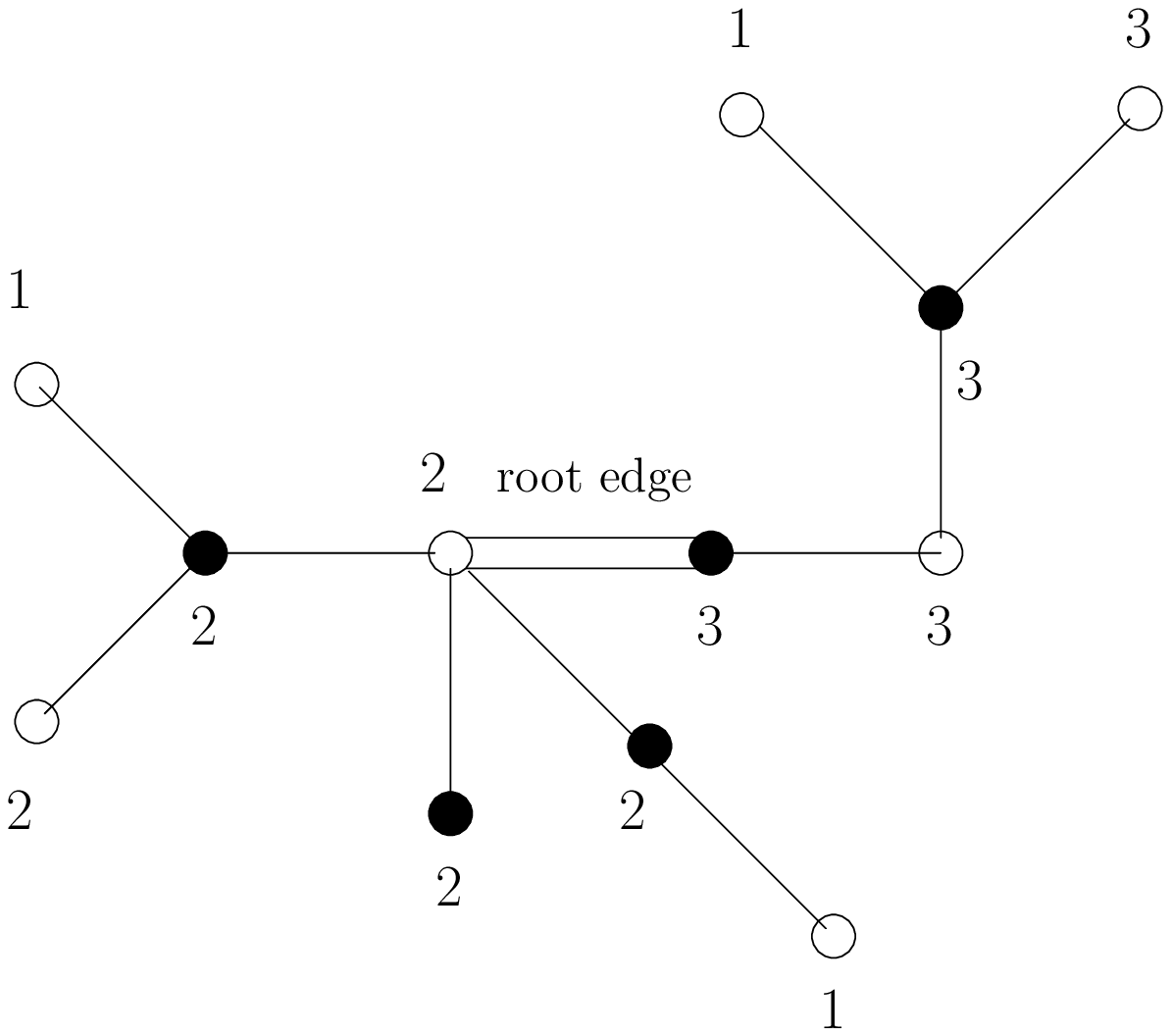}
  \caption{The coloured black and white plane edge rooted tree in this figure is the same as the one on the left in Figure \ref{fig:baseplanetree}, except with colours.  Its edges would be labelled as the tree on the right in Figure \ref{fig:baseplanetree}.  Thus, this tree corresponds to $(\alpha, \psi)$ and $(\beta, \nu)$, where $\alpha =
  (1\; 6\; 8\; 9)(2\; 5)(3)(4)(7)(10)(11)$ and $\beta = (1\; 5)(2\; 3\; 4)(6\;
  7)(8)(9\; 10\; 11)$, and where the cycles $(1\;
  6\; 8\; 9)$ and $(4)$ of $\alpha$ are coloured 2 and 3 by $\psi$,
  respectively.  Here, we assume $m \geq 3$.}\label{fig:colourplanetree}
\end{figure}

It is easy to see from the above construction, that top
\emph{coloured}
factorizations are obtained in the following way.  Let $(\alpha, \psi) \in
\colourS_k$ and $\beta$ be such that $\alpha \beta =
(1\; 2\; \cdots k)$.  Now, using the construction of Goulden and Jackson,
from $\alpha$ and $\beta$ create a black and white plane edge rooted tree.
As the white vertices of the tree correspond to cycles of $\alpha$, give the
white vertices an additional colour $i$ for $1 \leq i \leq m$ using
$\psi$.  For a black vertex $v$, an additional colour $i$ for $1 \leq
i \leq m$ is given with the rule that $v$ obtains colour $j$, where $j$ is the
maximum colour amongst all neighbours of $v$.  Thus, to be clear, vertices
have two types of colours;  they are either black or white and they have a
colour $i$, with $1 \leq i \leq m$.  As the black vertices
determine the cycles of $\beta$, the labels $i$ determine a function
$\phi: C_k(\beta) \longrightarrow [m]$.  Note that $\phi$ clearly determines a
function $\phi^\prime : C_k(\beta^{-1}) \longrightarrow [m]$ by
$\phi^\prime(u^{-1}) = \phi(u)$ for any cycle $u$ of $\beta$.  One can easily
check that $(\alpha, \psi) \circ (1\; 2\; \cdots k)^{-1} = (\beta^{-1},
\phi^\prime)$ (see Figure \ref{fig:colourplanetree}).  We will, therefore, call
the set of black and white plane edge rooted trees with this additional colour
restriction \emph{coloured black and white plane edge rooted trees} and
denote this class by \tclass.
Define 
for a tree $T \in \tclass$ the weights $\omega_w^{(m)}(T) =
(\omega_{w_1}^{(m)}(T), \omega_{w_2}^{(m)}(T),\dots)$ and $\omega_b^{(m)}(T) =
(\omega_{b_1}^{(m)}(T), \omega_{b_2}^{(m)}(T),\dots)$ where
$\omega_{w_i}^{(m)}(T)$ and $\omega_{b_i}^{(m)}(T)$ are the number of
white, respectively black, vertices in
$T$ coloured $i$.  As usual, let $\boldp^{\omega_w^{(m)}(T)} = \prod_i
p_i^{\omega_{w_i}^{(m)}(T)}$ and $\boldq^{\omega_b^{(m)}(T)} = \prod_i
q_i^{\omega_{b_i}^{(m)}(T)}$, and define
\begin{equation*}
	\gent = \sum_{T \in \tclass}
	\boldp^{\omega_{w}^{(m)}(T)} \boldq^{\omega_{b}^{(m)}(T)}
	x^{\textnormal{number of vertices of } T}.
\end{equation*}
Evidently, we have the following proposition from the above discussion.
\begin{proposition}\label{thm:deft}
\begin{equation*}
	\gent = \fact_{\boldpq}(x).
\end{equation*}
\end{proposition}
In the following sections we show that $\gent + (\sumofvars{p}{m})x - 1 = -
G_{\boldpq}(-x)$, which proves that \eqref{eq:gentopterms} holds by Proposition
\ref{thm:deft}.

\section{Planted Trees}

\begin{figure}[ht]
 \centering
  \includegraphics[width=9.0cm]{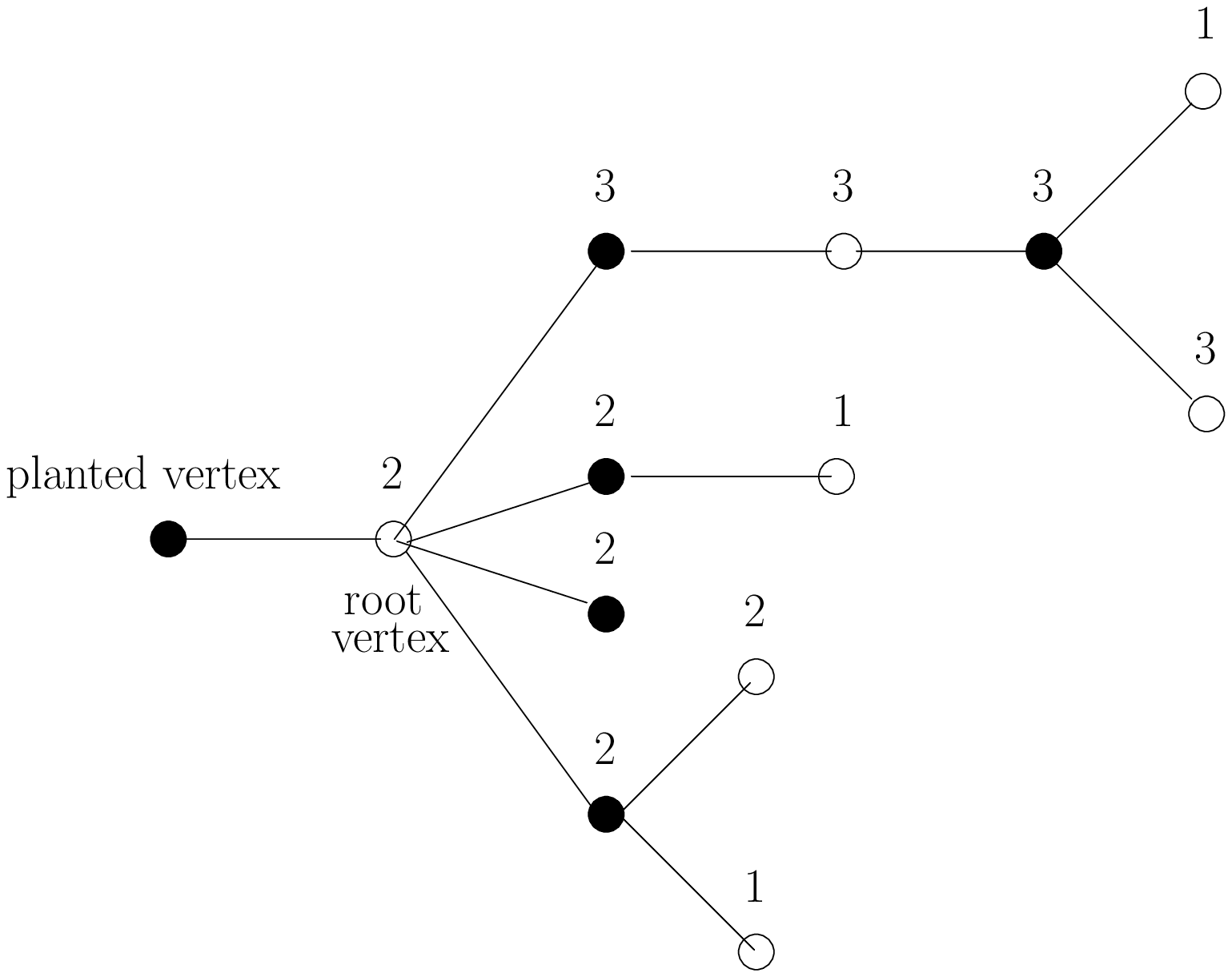}
  \caption{The coloured plane planted tree in $\bclass_2$ (here, we assume
  that $m \geq 3$) that corresponds to the coloured plane tree in Figure
  \ref{fig:colourplanetree}.  The tree is obtained by attaching a planted vertex
  to the white vertex incident with the root edge of the tree in Figure
  \ref{fig:colourplanetree}.  The root edge in Figure \ref{fig:colourplanetree}
  becomes the first of the linearly ordered edges (from top to bottom) emanating
  from the root.}\label{fig:colourplaneplantedtree}
\end{figure}

It is clear that the class of trees \tclass is in bijective correspondence with
the following class.  
Let $\bclass_i$ be the set of \emph{coloured} plane
\emph{planted} trees whose planted vertex is coloured black (the planted vertex does
not otherwise have a colour $i$), the vertex adjacent to the planted vertex,
which we call the \emph{root},
is white and coloured with the colour $i$ and the colouring of the rest
of the tree is consistent with the class of
trees in $\tclass$.  The planted vertex gives a linear order to the 
edges connecting the root to its children (see Figure \ref{fig:colourplaneplantedtree}).  Define $\wclass_i$ analogously;  that is,
$\wclass_i$ is the class of plane planted trees with planted vertex coloured white
(but with no colour $i$) and with the black root vertex coloured 
$i$.  In both these classes of trees, a vertex $v$ is the parent of a
vertex $w$ and $w$ is, likewise, called a child of $v$ if $v$ and $w$ are
connected by an edge and $v$ is on the unique path joining $w$ to the
planted vertex of the tree.  A tree in the class $\bclass_2$ is given in 
Figure \ref{fig:colourplaneplantedtree}.
Define the generating series
\begin{equation}
	\begin{split}
	B_i(x)&= \sum_{T \in \bclass_i} \boldp^{\omega_{w}^{(m)}(T)}
	\boldq^{\omega_{b}^{(m)}(T)} x^{\textnormal{number of non planted vertices of } T}\\
	W_i(x)&= \sum_{T \in \wclass_i} \boldp^{\omega_{w}^{(m)}(T)}
	\boldq^{\omega_{b}^{(m)}(T)} x^{\textnormal{number of non planted vertices of } T}
\end{split}
\end{equation}
If we let $\mathcal{P}_i$ be the class of white vertices with label $i$ then
it is easy to see that
\begin{equation*}
	\bigcup_{i=1}^m \mathcal{P}_i \bigcup \tclass = \bigcup_{i=1}^m \bclass_i,
\end{equation*}
from which it follows
\begin{equation*}
	\gent + (\sumofvars{p}{m})x = \sum_{i=1}^m B_i(x)
\end{equation*}
Thus, in order to find an expression for the generating series \gent we find one
for $\sum_{i=1}^m B_i(x)$ (see Figure \ref{fig:colourplaneplantedtree}).  For
convenience, we set $\geni$ be the previous generating series;  that is, we
set
\begin{equation}\label{eq:defi}
	\geni = \gent + (\sumofvars{p}{m})x = \sum_{i=1}^m B_i(x)
\end{equation}

We now find relations between the classes $\bclass_i$ and $\wclass_i$ for $1
\leq i \leq m$.  In order to do this we introduce a final class of trees.
Define the class of \emph{improperly} coloured trees, denoted $\wclasshat_i$,
planted at a white coloured vertex (but otherwise does not have a colour $i$)
and a black root coloured $i$.  The
white children of the black root vertex can
only be coloured with the colours of $1, 2, \dots, i-1$ (hence, the name
\emph{improperly} coloured).  Also note, we insist that the black root
has non-empty subtree below it (otherwise, such a tree would
\emph{not} be improperly labelled).  Define the generating series $\hat{W}_i(x)$
of the class $\wclasshat_i$ analogously to the series $W_i(x)$.  For $i=1$, the class $\wclasshat_i$ is empty
and its corresponding generating series is $\what_1(x) = 0$.  We shall see their importance in the next section.

\section{Decomposition of the classes $\bclass_i$ and $\wclass_i$ and the
proof of the Main Theorems}\label{sec:proofmain}

\begin{figure}[ht]
 \centering
  \includegraphics[width=9.0cm]{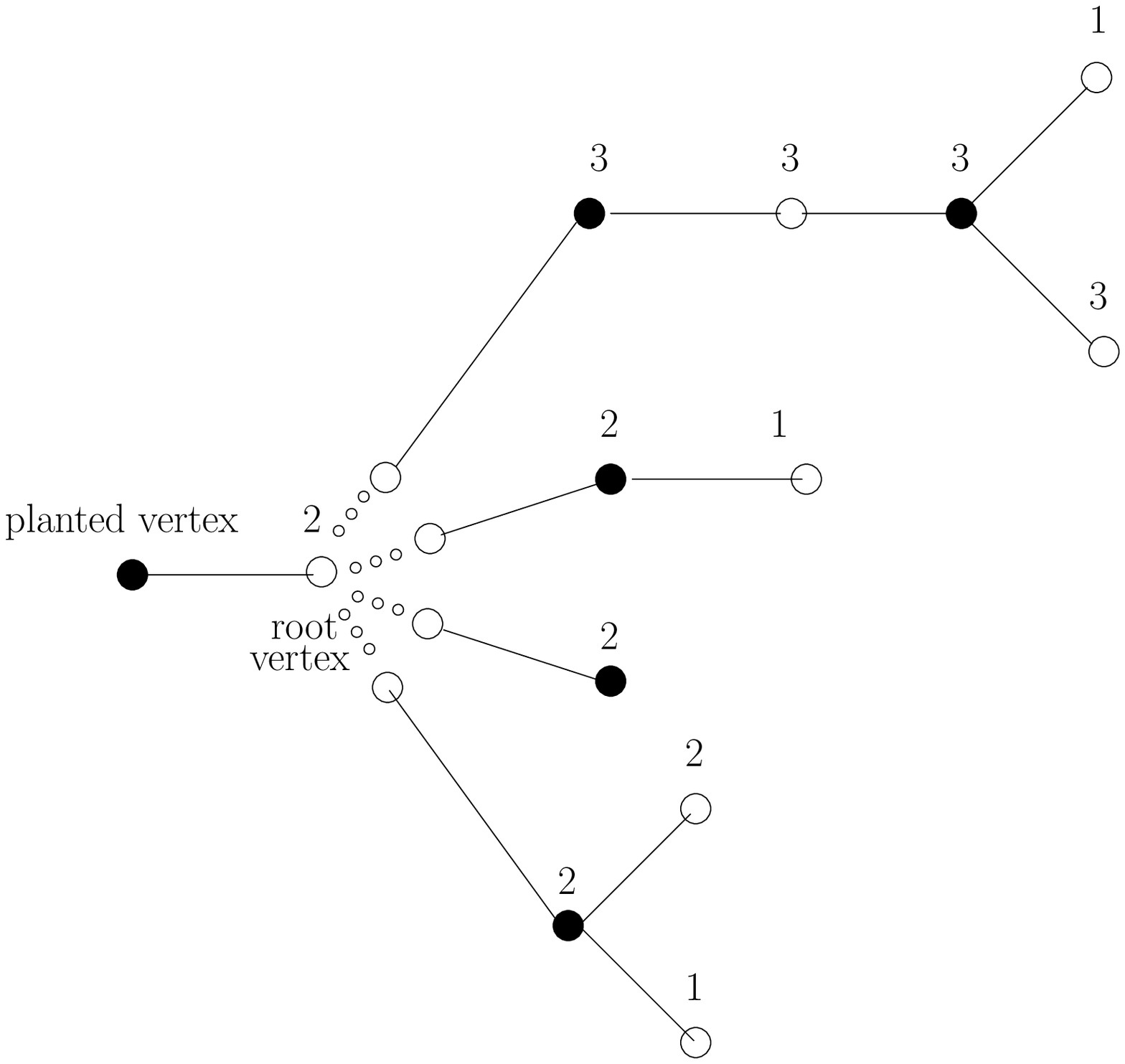}
  \caption{The coloured plane planted tree in $\bclass_2$ given in Figure
  \ref{fig:colourplaneplantedtree} decomposed into two $\wclass_2$, one $\wclass_3$
  and one $\wclasshat_2$ trees.}\label{fig:colourplaneplantedtreebroken}
\end{figure}

We begin by discussing the decomposition of the class $\bclass_i$, for $1 \leq i \leq m$.  Recall, a tree
in this class has a planted black vertex adjacent to a white root vertex with
the colour $i$.   Every child of the root vertex is
black and because of the colouring rule requiring the colour of a black vertex
to be the largest colour amongst it white neighbours, the colours $i, i+1, \dots,
m$ are the possible colours for the children of the white root vertex.  Each of these black children have subtrees to which they are
attached, and can be made into a planted $\wclass_j$ tree for $i \leq j \leq m$
by attaching a planted white vertex to each of the black children (see, for
example, Figure \ref{fig:colourplaneplantedtreebroken}).
Note, however, we may also obtain an improperly labelled tree $\wclasshat_j$ for
$i \leq j \leq m$ (as in, for example, the second subtree counting from the top
in Figure \ref{fig:colourplaneplantedtreebroken}).  The final caveat here is
that any black child of the white root vertex with colour strictly greater than $i$, must have a non-trivial subtree below it (for
otherwise that black vertex would be not be properly coloured).
Since the subtrees are linearly ordered, we see that $B_i(x)$ satisfies
\begin{equation}\label{eq:decompbi}
	B_i(x) = \frac{p_i x}{1 - \left(W_i(x) + \what_i(x) + \left(W_{i+1}(x) -
	q_{i+1} x \right) + \cdots + \left(W_m(x) - q_m x\right) \right)}.
\end{equation}
We repeat that $\what_1(x) = 0$.
We can similarly find expressions for the generating series $\what_i(x)$ and
$W_i(x)$.  Beginning with $\what_i(x)$, for $2 \leq i \leq m$, we have by definition a tree in the class 
$\wclasshat_i$ has a white planted
vertex and black root vertex with colour $i$.  The black root
vertex has white children coloured with $1, 2, \ldots i-1$.  Also by
definition, the black root vertex must have a non-trivial subtree beneath it.
Since the subtrees are linearly ordered, we have
\begin{equation}\label{eq:defwhat}
	\what_i(x) = \frac{(q_i x) (B_1(x) + B_2(x) + \cdots + B_{i-1}(x))}{1 -
	(B_1(x) + B_2(x) + \cdots + B_{i-1}(x))}.
\end{equation}
For the trees $\wclass_i$, by definition they are trees with a white planted
vertex and a black root vertex with colour $i$ and are properly coloured;  that
is, the black root vertex has no white children with colour greater
than $i$.  The black root vertex may have a trivial subtree
beneath it, but if it does not it must have at least one white child coloured
$i$ (in order be properly coloured).  As these are the only restrictions, we see
\begin{align}
	W_i(x) &= \frac{q_i x}{1 - (B_1(x) + B_2(x) + \cdots + B_i(x))} -
	\frac{(q_i x) (B_1(x) + B_2(x) + \cdots + B_{i-1}(x))}{1 - (B_1(x) +
	B_2(x) + \cdots + B_{i-1}(x))}\notag\\
	&= \frac{q_i x}{1 - (B_1(x) + B_2(x) + \cdots + B_i(x))} -
	\what_i(x).\label{eq:defw}
\end{align}
We now show that $\geni$ given in \eqref{eq:defi} satisfies the same equation as
\eqref{eq:stantopgoodform}.

\begin{lemma}
	\label{thm:firstrec}
	\begin{equation*}
		B_i(x) = \frac{p_i x (B_1(x) + \cdots + B_i(x)
		-1)}{\geni - 1 - (p_{i+1} + \cdots + p_m)x + q_i x}.
	\end{equation*}
\end{lemma}
\begin{proof}
	Our proof is by induction on $i$ beginning at $i=m$.
	For $i=m$ we have by \eqref{eq:decompbi}
	\begin{align*}
		B_m(x) &= \frac{p_m x}{1 - (W_m(x) + \what_m(x))}\\
		&= \frac{p_m x}{1 - \frac{q_m x}{1 - (B_1(x) + \cdots +
		B_m(x))}}\\
		&= \frac{p_m x (B_1(x) + \cdots + B_m(x) - 1)}{\geni - 1 + q_m
		x},
	\end{align*}
	completing the base case.

	Now suppose that our conclusion is true for $i=t$.  We wish to show our
	conclusion holds for $i= t-1$.  By the induction hypothesis, we have
	\begin{equation*}
		B_t(x) = \frac{p_t x (B_1(x) + \cdots + B_t(x)
		-1)}{\geni - 1 - (p_{t+1} + \cdots + p_m)x + q_t x}.
	\end{equation*}
	and from \eqref{eq:decompbi} we have
	\begin{equation*}
		B_{t-1}(x) = \frac{p_{t-1} x}{1 - \left(W_{t-1}(x) +
		\what_{t-1}(x) +
		\left(W_{t}(x) - q_{t} x \right) + \cdots + \left(W_m(x) -
		q_m x\right) \right)}.
	\end{equation*}
	from which we obtain
	\begin{align}
		B_{t-1}(x) &= \frac{-p_{t-1} x}{\left(W_{t-1}(x)
		+ \what_{t-1}(x) - \what_t(x) - q_t x
		\hspace{5.0cm}\right.}\notag\\
		&\hspace{2.0cm}\left.
		+\left(W_t(x) + \what_t(x) + W_{t+1}(x) - q_{t+1} x +
		\cdots + W_m(x) - q_m x \right)\right) -1\notag\\
		&=\frac{-p_{t-1} x}{\frac{q_{t-1} x}{1 - \left(B_1(x) + \cdots +
		B_{t-1}(x) \right)} - \frac{q_t x}{1 - \left(B_1(x) + \cdots +
		B_{t-1}(x) \right)} + \frac{\geni -1 - \left(p_{t+1} + \cdots +
		p_m\right) x + q_t x}{1 - \left(B_1(x) + \cdots +
		B_{t}(x) \right)}},\label{eq:mess}
	\end{align}
	where the last summand in the denominator of \eqref{eq:mess} follows
	from \eqref{eq:decompbi} with $i=t$ and the induction hypothesis.  Continuing to simplify,
	\eqref{eq:mess} becomes
	\begin{align*}
		&\frac{-p_{t-1} x(1 - (B_1(x) + \cdots +B_t(x)))(1 - (B_1(x) +
		\cdots +B_{t-1}(x)))}{(1 - (B_1(x) + \cdots +
		B_t(x)))\left((q_{t-1} x - q_t x) + \left(\geni - 1
		\right.\right.\hspace{4.4cm}}\\
		&\hspace{1.0cm}\left.\left. - (p_{t+1} + \cdots + p_{m})x + q_t
		x\right)\right) + B_t(x)(\geni - 1 - (p_{t+1} + \cdots + p_{m})x
		+ q_t x)\\
		&= \frac{-p_{t-1} x(1 - (B_1(x) + \cdots +B_t(x)))(1 - (B_1(x) +
		\cdots +B_{t-1}(x)))}{(1 - (B_1(x) + \cdots +
		B_t(x)))(\geni - 1 -
		(p_{t+1} + \cdots + p_{m})x + q_{t-1} x)\hspace{2.0cm}}\\
		&\hspace{7.0cm}+p_t x(B_1(x) + \cdots + B_t(x) -1)\\
		&= \frac{p_{t-1} x(B_1(x) + \cdots + B_{t-1}(x) -1) }{\geni - 1 -
		(p_t + p_{t+1} + \cdots + p_m)x + q_{t-1}x},
	\end{align*}
	completing the proof.
\end{proof}
\begin{lemma}
	For $0 \leq i \leq m$ the following equation holds:
	\begin{equation*}
		B_1(x) + \cdots + B_i(x) -1 = (\geni -1)\prod_{j=i+1}^m \frac{
		\left(\geni -1 -(p_j + \cdots +p_m)x +
		q_jx\right)}{\left(\geni -1 - (p_{j+1} + \cdots
		+ p_m)x + q_jx\right)}
	\end{equation*}
	\label{thm:seclem}
\end{lemma}
\begin{proof}
	The proof is by induction on $i$, beginning at $i=m$.  The case $i=m$ is
	trivial.

	Now suppose for $i=t$ the statement of this lemma is true;  that is,
	\begin{equation*}
		B_1(x) + \cdots + B_t(x) -1 = (\geni -1)\prod_{j=t+1}^m \frac{
		\left(\geni -1 -(p_j + \cdots +p_m)x +
		q_jx\right)}{\left(\geni -1 - (p_{j+1} + \cdots
		+ p_m)x + q_jx\right)}.
	\end{equation*}
	By Lemma \ref{thm:firstrec}, we have
	\begin{equation*}
		B_t(x) = \frac{p_t x \left(B_1(x) + \cdots + B_t(x)
		- 1\right)}{\geni -1 - (p_{t+1} + \cdots + p_m)x + q_tx)}.
	\end{equation*}
	But,
	\begin{align*}
		B_1(x) + \cdots &+ B_{t-1}(x) -1\\
		&= B_1(x) + \cdots + B_{t}(x)
		-1 - \frac{p_t x \left(B_1(x) + \cdots + B_t(x) -1
		\right)}{\geni -1 - (p_{t+1} + \cdots + p_m)x + q_tx}\\
		&= \left(B_1(x) + \cdots + B_{t}(x) -1 \right) \left( 1 - \frac{p_t x}{\geni
		-1 - (p_{t+1} + \cdots + p_m)x + q_tx}\right)\\
		&= \left( (\geni -1)\prod_{j=t+1}^m \frac{
		\left(\geni -1 -(p_j + \cdots +p_m)x +
		q_jx\right)}{\left(\geni -1 - (p_{j+1} + \cdots
		+ p_m)x + q_jx\right)}\right)\\
		& \hspace{4.0cm} \cdot \frac{(\geni -1 - (p_t + \cdots + p_m)x +
		q_t x)}{(\geni -1 -(p_{t+1} + \cdots + p_m)x + q_tx)}\\
		& = (\geni -1)\prod_{j=t}^m \frac{
		\left(\geni -1 -(p_j + \cdots +p_m)x +
		q_jx\right)}{\left(\geni -1 - (p_{j+1} + \cdots
		+ p_m)x + q_jx\right)},
	\end{align*}
	completing the proof.
\end{proof}
We now give a proof of the main theorems.
\begin{proofof}[Proof of Theorem \ref{thm:main1}.]
	From Lemma \ref{thm:seclem}, we have by setting $i=0$
	\begin{equation*}
		(\geni -1)\prod_{j=1}^m
		\frac{ \left(\geni -1 -(p_j + \cdots +p_m)x +
		q_jx\right)}{\left(\geni -1 - (p_{j+1} + \cdots + p_m)x +
		q_jx\right)} = -1.
	\end{equation*}
	Thus, $\geni - 1$ and $-G_{\boldpminq}(-x)$ both satisfy
	\eqref{eq:stantopgoodform}.  Combining this with \eqref{eq:defi}
	and Proposition \ref{thm:deft} gives the result.
\end{proofof}
\begin{proofof}[Proof of Corollary \ref{thm:main2}]
If $\mu = \mu_1 \mu_2 \ldots \mu_\ell \vdash k$, then the terms of highest degree in
$(-1)^k F_\mu(\boldpminq)$, which have degree $k+ \ell(\mu)$,
are given by $\prod_{i=1}^\ell (-1)^{\mu_i} G_{\mu_i}(\boldpminq)$ (see \citet[Theorem
4.9]{sniasymp} and \cite[Theorem 9]{snialarge} and references therein.  To see how Kerov polynomials are used to obtain characters of the
symmetric group and Stanley's polynomials, see \citet{ratt:6}).  
Assuming that $\alpha \beta = \gamma$ in $\mfS_k$ where $\gamma$ has cycle type $\mu$
and $\kappa(\alpha) + \kappa(\beta) = k + \kappa(\gamma)$ then the product $\alpha
\beta = \gamma$ necessarily decomposes into $\kappa(\gamma)$ products of the form
$\alpha_i \beta_i = \gamma_i$, where 
\begin{enumerate}
	\item $\gamma_i$ is a cycle of $\gamma$;
	\item $\alpha_i = \alpha_i^1 \alpha_i^2 \cdots \alpha_i^s$ and $\beta_i = \beta_i^1 \beta_i^2 \cdots \beta_i^t$, where $\alpha_i^j$ and $\beta_i^j$ are cycles in $\alpha$ and $\beta$,  respectively; the cycles $\alpha_i^j$ and $\beta_i^j$ are precisely the cycles in $\alpha$ and $\beta$ that contain elements in the support of $\gamma_i$;
i.e. the elements $h$ of $\{1, 2, \ldots, k\}$ such that $\gamma_i(h) \neq h$;
	\item if the length of $\gamma_i$ is $k_i$ then $\kappa(\alpha_i) +
		\kappa(\beta_i) = k_i + 1$
\end{enumerate}
(see
\citet[Section 4]{goul:7}).   From this it follows that
\begin{equation*}
	\sum_{(\alpha, \psi) \in \colourS_k \atop \kappa(\alpha) +
	\kappa(\alpha \omega_\mu) = k+\kappa(\mu)}
		\boldp^{\km(\alpha, \psi)} \boldq^{\km( (\alpha, \psi)\circ
		\omega_\mu)} = \prod_{i=1}^{\ell(\mu)}(-1)^{\mu_i}
		G_{\mu_i}(\boldpminq).
\end{equation*}
completing the proof.
\end{proofof}

\begin{center}
	{\bf \large{Acknowledgements}}
\end{center}
	This work was supported by a \emph{Natural Sciences and Engineering
	Research Council of Canada} Postdoctoral Fellowship.  I would like
	to thank Richard Stanley for communicating his conjecture to me and for
	some helpful discussions.  I would also like to thank Karola Meszaros
	and the anonymous referee for their useful comments on the previous
	version of this manuscript.

\bibliographystyle{plainnat}
\bibliography{colourfacts}
\end{document}